\documentclass[preprint]{elsarticle}
\usepackage{amsmath}
\usepackage{amssymb}
\usepackage{amsthm}
\usepackage{amsfonts}
\usepackage{amsbsy}
\usepackage{amstext}
\usepackage{amsfonts}
\usepackage{amstext}
\usepackage{mathrsfs}
\usepackage{graphpap}
\usepackage{graphics,graphicx}
\usepackage[T1]{fontenc}
\usepackage{subcaption}
\usepackage{multicol}
\usepackage{booktabs}
\usepackage{array}
\usepackage{bm}
\usepackage{siunitx}
\usepackage{hyperref}
\usepackage[section]{placeins}
\usepackage{xargs}
\usepackage{algorithmic,algorithm}

\usepackage{dirtytalk}
\usepackage[pdftex,dvipsnames]{xcolor}
\usepackage[colorinlistoftodos,
   prependcaption,
   textsize=tiny]{todonotes}
\usepackage{mathtools, cuted}
\usepackage{hyperref}
\usepackage[nameinlink,noabbrev]{cleveref}
\usepackage{lineno}

\newcommand{\IR}{{\mathbb R}}
\newcommand{\R}{{\mathbb R}}
\newcommand{\IN}{{\mathbb N}}

\newcommand{\beq}{\begin{equation}}
\newcommand{\eeq}{\end{equation}}
\newcommand{\bal}{\begin{align}}
\newcommand{\eal}{\end{align}}
\newcommand{\beqn}{\begin{equation*}}
\newcommand{\eeqn}{\end{equation*}}
\newcommand{\baln}{\begin{align*}}
\newcommand{\ealn}{\end{align*}}

\newcommand{\Pb}{\mathbb P}

\newcommand{\E}{\mathbb{E}}

\newcommand{\I}{{\mathbb I}}
\usepackage{xcolor}
\newtheorem{remark}{Remark}
\newtheorem{definition}{Definition}
\newdefinition{lemma}{Lemma}
\hypersetup{colorlinks=false,allbordercolors=blue,pdfborderstyle={/S/U/W 1}}

\begin{document}
\markboth{F. Baltazar-Larios, F. Delgado-Vences and G. Salcedo-Varela}{Simulating diffusion bridges using the Wiener chaos expansion}

\begin{frontmatter}
    \title{Simulating diffusion bridges using the Wiener chaos expansion}
        \author[1]{Francisco Delgado-Vences}
        \ead{delgado@im.unam.mx}
        \author[2]{Gabriel Adri\'an Salcedo-Varela}
        \ead{avarela55@gmail.com}
        \author[3]{Fernando Baltazar-Larios}
        \ead{fernandobaltazar@ciencias.unam.mx}
    \affiliation[1]{organization={CONAHCYT-UNAM, Instituto de Matem\'aticas}, 
    addressline={Sede Oaxaca}, city={Oaxaca de Ju\'arez}, country={M\'exico}}
    \affiliation[2]{organization={CONAHCYT-UNAM, Facultad de ciencias}, 
    addressline={Ciudad Universitaria}, city={Ciudad de M\'exico}}
    \affiliation[3]{organization={Facultad de ciencias-UNAM}, 
    adressline={Ciudad Universitaria},city={Ciudad de M\'exico}}
    \begin{abstract}
In this paper, we present a simulation method for diffusion bridges by utilizing an approximation of the Wiener-Chaos Expansion (WCE) applied to a related diffusion process. Indeed, we consider the solution of a proposed diffusion bridge defined through a stochastic differential equation (SDE), and apply the WCE to this solution. As a result, we have developed a rapid method for simulating proposal diffusion bridges that eliminates the need for rejection rates. The method presented in this work could be beneficial in statistical inference, where many parameters are unknown and do not have associated information. This initial uncertainty can influence the effectiveness of data augmentation methods used in inference and increase their rejection rates. In opposition, our method excels by reliably constructing diffusion bridges quickly and efficiently, regardless of initial parameter conditions. We validate our method with a simple Ornstein-Uhlenbeck process. We have also compared our method with two existing methods and determined that its performance is satisfactory. We apply our method to four examples of SDEs and show the numerical results.
    \end{abstract}
\end{frontmatter}

    \section{Motivation and Introduction}  

    Stochastic differential equations (SDEs) play a crucial and flexible role in modeling dynamic phenomena, capturing the influence of random noise on temporal evolution.  Indeed, SDEs effectively capture both the underlying dynamics and the inherent randomness of actual observations; they have been rigorously applied to a diverse range of phenomena characterized by noise, such as finance, atmospheric science, biology, physics, ecology, genetics, oceanography, psychology, and medicine(see \cite{mao}, \cite{oksendal2013stochastic}, \cite{panik2017stochastic} for examples).
    
Let $X(t)$ represent the state of the system at time $t\geq 0$, which allows us to characterize its evolution over time using the following SDE
\begin{equation}\label{SDE1}
         \begin{aligned}
            dX(t)
                & = 
                    f(X(t))dt + \sigma(X(t))dB(t), 
         \end{aligned}
     \end{equation}

where $f$ (drift coefficient) represents the function of change at time $t$, $\sigma$ (diffusion coefficient)  is a known function,  and $\{B_t\}_{0\leq t\leq T}$ is a standard Brownian process which is interpreted as the input random noise in the system. It is known that $f$ models the mean of the process, or the deterministic trend of the data (see \cite{oksendal2013stochastic}, \cite{panik2017stochastic}). Furthermore, it is commonly assumed that $f$ and $\sigma$ include certain parameters, denoted as $\Theta$, which are used to fit the model \eqref{SDE1} to actual data. Consequently, it becomes imperative to estimate these parameters $\Theta$ using available observations. Estimating parameters generally requires the application of statistical inference techniques.

  In statistical inference, it is generally assumed that we have sufficient information/data about a phenomenon to effectively employ a variety of statistical methods. Nonetheless, in many cases, the observations may be insufficient for effectively applying statistical methods, which can lead to inaccurate inference results. To address this challenge, an effective strategy entails filling in the gaps between the information provided; that is, data augmentation. In particular, a popular approach, in the case  of parameter inference for discretely observed stochastic differential equations (SDEs), involves the use of {\it diffusion bridges} to complete the information. This methodology successfully creates a stochastic path by linking all actual observations through the simulation of a conditioned diffusion process.


In this work, we address the problem of simulating a one-dimensional diffusion process $\bm{X}=\{X(t)\}_{t=0}^ T$ given by the SDE \eqref{SDE1}  subject to the initial and final conditions
    \begin{equation}\label{SDE_conditions}
    X(0) =\eta,\quad
    X(T)= \theta.
    \end{equation} We assume that $f$ and $\sigma$ satisfy suitable conditions to ensure the 
    existence and uniqueness of the solution $\bm{X}$ in some sense (strong or weak).

   \paragraph{Literature review} 
   
    The simulation of diffusion bridges represents a considerable challenge, as the pursuit of efficient methodologies to accurately simulate these processes has proven to be exceedingly demanding. This area of research has undergone active investigation over the past two decades. Diffusion bridge simulation is used to complete missing paths between observations at discrete times. One of the most used methodologies is the Expectation-Maximization algorithm or Markov Chain Monte Carlo method (MCMC) to solve the problem of statistical inference for the SDE  \eqref{SDE1} when missing data. 
    
   The list of references provided herein is far from exhaustive; we direct the reader to \cite{Val-21} and \cite{chau2024efficient} for a more comprehensive compilation of references. To our knowledge, the first methods for diffusion bridge simulation were presented in \cite{clark1990simulation}, who introduce the proposal bridge that is similar to the one that we are studying in this paper. Later on, \cite{durh-02} and \cite{Roberts-01}, present works that are based on Metropolis–Hastings algorithm. Algorithms for exact simulation of diffusion bridges were developed in \cite{Beskos08} and \cite{Beskos06}, their methods can be employed to sample diffusion bridges without any time-discretization error but are limited to the class of diffusion processes that can be transformed into one with a unit diffusion coefficient. 
    
    In \cite{Lin-10}, the authors proposed a method based on Monte Carlo simulations using the empirical distribution of backward paths. Related works inspired by the same idea include guided proposal densities (\cite{Del-15} and \cite{Schauer17}), sequential Monte Carlo algorithms incorporating resampling with backward information filter approximations (\cite{Gua-17});  the solution of the backward Kolmogorov partial differential equation (\cite {Wan-20}); and a method based on coupling and time-reversal of the diffusion (\cite{bladtSoren} and \cite{bladtSorenc});  This last method is simple to implement and demonstrates good performance over extended time periods. However, it requires the existence of the invariant distribution and the ability to sample from it. In \cite{Val-21}, the authors show that the time-reversed diffusion bridge process can be simulated if one can time-reverse the unconditioned diffusion process. They also present a variational formulation to learn this time-reversal that relies on a score-matching method to circumvent intractability. In \cite{chau2024efficient}, the authors use an approximation for the transition density of the solution of the related SDE $X$.  Their approximation  relies on the Itô-Taylor expansion applied to a smooth approximation of the Dirac function, which is applied to
  the infinitesimal generator of $X$. Subsequently, they employ the Euler/Milstein scheme to simulate the diffusion bridge.

Recently, the authors of \cite{bierkens2021piecewise} proposed a method based on an approximation for the density of the SDE's solution. 
They use as a main tool a truncation of the Faber–Schauder basis, which is similar to the Lévy–Ciesielski construction of a Brownian motion. The underlying intuition is that, given mild assumptions on the drift function, any diffusion process exhibits local behavior similar to that of a Brownian motion. Their proposed method may be regarded as a spectral method, as once the basis is established—the Faber-Schauder basis in that instance—the only task that remains is to determine the coefficients, which are random in that particular case, to settle the approximation. In contrast to the previous work, our proposed method utilizes a stochastic basis, then we fix the coefficients that are time-dependent functions, closely mirroring the dynamics of diffusion.

 \paragraph{Our approach}    In this paper, we use an approximation of the Wiener-Chaos Expansion (WCE) to propose a novel method to simulate diffusion bridges. The WCE is a spectral decomposition in the random parameter; meaning, that if a random variable is square integrable, thus we could write it as a Fourier series, with the stochastic basis being a family of Hermite polynomial functional evaluated in a sequence of independent Gaussian random variable (see below for further details). Furthermore, in our specific case, the coefficients constitute a sequence of deterministic functions that must be defined to carry out the spectral decomposition. Thus, the idea of this work is to consider the solution of SDEs and find the WCE of a particular representation of the diffusion bridge. An advantage of the method presented here is that the numerical performance to simulate the diffusion bridges works well for any time interval between observations (long of the bridge). Another advantage is that the execution time for simulating a diffusion bridge is the same regardless of whether the estimators used are close, or not, to the true parameter values. This fact is highly significant in statistical inference, particularly when simulating diffusion bridges with estimators that may deviate significantly from the theoretical true parameter values.
    
    To propose our method,  let us denote by $(\Omega, \mathcal{F}, \Pb)$  the canonical probability space 
    associated with a one-dimensional standard Brownian motion $\{B(t)\}_{t=0}^T$. This means $\Omega =C_0([0,T]; \IR)$, $\Pb$ is the Wiener measure, and $\mathcal{F}$ is the completion of the Borel $\sigma$-field of $\Omega$ 
    with respect to $\Pb$. Set the Hilbert space $\mathbb{H}=L^2([0,T])$. In particular, we suppose that $\bm X\in L^2(\Omega; \mathbb{H})$. This assumption 
    allows us to use the spectral representation of the solution using the WCE (cf. \cite{lototsky2006stochastic} or 
    \cite{lototsky2017stochastic}).
    
    Consider the integral version of the SDE \eqref{SDE1}:
    \begin{equation}\label{integral_SDE}
            X(t) = x_0+ \int_0^t f(X(s))ds +  \int_0^t  \sigma(X(s))dB_s. 
    \end{equation}
    Observe that in \eqref{integral_SDE} the initial condition is added naturally, 
    but the final condition $X(T)=x_T$ is not. 
    
    The type of stochastic boundary problem given by \eqref{SDE1}-\eqref{SDE_conditions} 
    has been studied before and is a challenging problem. However, in some cases, 
    the existence and uniqueness of the solution have been shown by using sophisticated tools from stochastic analysis
    (cf \cite{nualart1991boundary}). In this work, we use 
    a diffusion $\bm Y=\{Y(t)\}_{t=0}^T$ (which can be seen as a proposal of the diffusion bridge) 
    such that when it is applied the Radon–Nikodym derivative to the laws of $\bm X$ 
    and $\bm Y$, the choice of drift function allows us to have that the law of $\bm Y$ is 
    absolutely continuous respect to  the law of $\bm X$ (see \cite{delyon2006simulation} for the hypothesis of $f$ and $\sigma$ that ensure this fact). The method studied here avoids the use of the Skorohod integral or enlargement of filtrations to make measurable the final condition.  
    
    Consider the diffusion  $\bm Y$ given by (see \cite{delyon2006simulation} for a general diffusion or \cite{Gasbarra2007}  for the Gaussian case )
    \begin{equation}
        \begin{aligned}\label{Yt-eq}
            dY(t) & =     (\theta-Y(t)) \frac{1}{T-t} dt + dX(t),\\
            Y(0) &=\eta.
        \end{aligned}
    \end{equation}
     Observe that since a final 
    value for the diffusion is included for the SDE $\bm X$, $X(T)=\theta$, then the 
    diffusion $\bm X$ is an anticipative diffusion or non-adapted solution. This implies 
    the necessity to use the Skorohod integral to define properly the solution of 
    the SDE \eqref{SDE1} (see \cite{nualart1991boundary}) or to enlarge the filtration. However, using $\bm Y$ we avoid these issues.
     By assuming the same conditions on $f$ and $\sigma$ as for $\bm X$, we have that $\bm Y\in L^2(\Omega; \mathbb{H})$, which allows us to use the spectral representation of the solution through the WCE for $\bm Y$.
     The key idea of this paper is to truncate the WCE for the diffusion $\bm Y$ and use it to approximate the diffusion bridge $\bm X$ given by the equations \eqref{SDE1}-\eqref{SDE_conditions}. 

 \paragraph{Our contribution}
In this paper, we consider a proposal bridge diffusion, which is absolutely continuous with respect to a bridge diffusion, and we apply the WCE to approximate it; indeed, the proposed numerical method is a truncation of the WCE for a particular SDE, that allows us to perform simulations of diffusion bridges in a highly efficient manner. In order to implement the proposed simulation method,a set of coupled ODEs must be solved. This system of equations can be solved iteratively; however, numerical methods are often required to obtain their solutions.\\

    The method takes advantage of the fact that the first ODE is actually the mean of the diffusion bridge. The remaining ODEs system constructs a deterministic bridge between $(0,0)$ to $(0,T)$. Consequently, when these equations are multiplied by the noises, they introduce a stochastic perturbation to the mean of the bridge. Therefore, we have a method that always builds a bridge between the endpoints, regardless of the distance or time jumps. In particular, our method does not utilize the existence of rejection probabilities in the construction of the diffusion bridges.\\

    The rest of the paper is organized as follows. In Section \ref{Malliavin-section} we review some facts from Malliavin calculus and the Wiener chaos expansion. Section \ref{WienerChaos-BM} applies the WCE to a particular SDE that will allow us to simulate diffusion bridges. We also present a truncation of the WCE and an error estimation for the aforementioned truncation. In Section \ref{validation_WCE}, we validate the proposed method; to do that we consider the Ornstein-Uhlenbeck process; and we construct some bridges between several points and we compare our method with the exact simulation method.  Section \ref{numerics}, contains the comparison of our method with other two numerical methods in the literature. Indeed, we compare two methods for the Ornstein-Uhlenbeck process and the Geometric Brownian motion. We construct some bridges between several points, and conduct a comparative analysis of our methodology against Bladt and S\o rensen's approach, and Doob's $h$-transform approach. Section \ref{numerics2}, contains the results of the numerical experiments when we apply our method to two SDEs. In Section \ref{conclusions} we present some concluding remarks of the method and extensions of the ideas introduced in this work.
    
    \section{A review on Malliavin calculus} \label{Malliavin-section}
    For the sake of completeness, we review some topics of Malliavin calculus 
    applied to our particular case. We also provide a brief review of the 
    WCE, and then we will apply it to the system of SDEs that 
    we have presented before.
    
    The WCE is a spectral decomposition in the probability space,
    which was introduced by Cameron and Martin in 1947, sometimes it is called the Fourier-Hermite expansion. Later on, it was revisited 
    by many researchers (see for instance \cite{luo2006wiener} for a good list of 
    references). The construction presented here was originally developed by Rozovzky 
    and Lototsky and their collaborators. It has been applied to SPDEs for theoretical 
    results such as existence and uniqueness, but also to provide a nice and efficient 
    framework to simulate solutions of SPDEs. For a theoretical introduction to SPDEs 
    and their solution via the WCE, we refer to 
    \cite{lototsky2017stochastic} and the references therein. For numerical purposes, we refer to the book \cite{zhang2017numerical} and the references therein. For the particular version of the WCE used in this manuscript, we refer to 
    \cite{lototsky2006stochastic}. For Malliavin calculus, see the classical 
    reference \cite{nualart2006malliavin}.
    
    Set $\mathbb{H} =L^2 ([0, T ])$ be the separable Hilbert space with inner 
    product $\langle \cdot, \cdot \rangle_\mathbb{H} $ and norm 
    $\|\cdot \|_\mathbb{H} $\footnote{Recall that the norm in this space is 
    given by $\|f \|_\mathbb{H}^2 = \int_0^T f^2(t) dt $.}. Denote by $\{e_i \}_{i\ge 1}$ an orthonormal basis of $\mathbb{H}$. We consider the isonormal process 
    $\bm W = \{W (h) : h \in \mathbb{H} \}$ indexed by $\mathbb{H}$
    defined on a probability space $(\Omega, \mathcal{F},\Pb)$, i.e. the random 
    variables $W (h)$ are centered Gaussian with a covariance structure 
        \begin{equation*}
            \E(W(h)W(g))= \langle h,g \rangle_\mathbb{H}, \quad 
            \mbox{for all } h,g\in\mathbb{H}.
        \end{equation*}
    We will assume that $\mathcal{F}=\sigma(\bm W)$. The key example of an isonormal process
    is the Wiener-It\^o integral $W(h)= \int_0^T h(s) dB(s)$, where $B(s)$ is the Wiener
    process and $h\in \mathbb{H}$. In fact, this will be the definition of isonormal process
    in our setting.
    
    Denote by $\{w_k(t)\}_{k\ge 1}$   a system of one-dimensional independent standard 
    Wiener processes.
    %
    %
    
    We now focus on the space $L^2(\Omega)=L^2(\Omega, \mathcal{F},\Pb)$ and we 
    will define an orthonormal basis for this space. First, we recall the real 
    normalized Hermite polynomials 
    \begin{equation*}
        H_0(x):=1, \qquad H_n(x):= \frac{(-1)^n}{\sqrt{n!}}
        \exp\Big(\tfrac{x^2}{2}\Big)\frac{d^n}{dx^n}
        \left(
            \exp\Big(-\tfrac{x^2}{2}\Big)
        \right).
    \end{equation*}
    Denote by $\mathcal{J}$ the set of multi-indices $m = (m_{k}: {k\ge 1} )$ 
    of finite length  $|m| := \sum_{k= 1}^{\infty} m_{k}<\infty$, i.e.
    \begin{equation*}
        \mathcal{J} = \Big\{ m = (m_{k})_{k\ge 1},\, m_{k}\in\{0, 1, 2,\ldots \}, 
        \,|m| < \infty\Big\}.
    \end{equation*}
    $\mathcal{J}$ is formed for a countable number of infinite dimensional vectors
    such that if $m\in \mathcal{J}$, then only a finitely many $m_{k}$ are not equal
    to zero.\\ 
    For $m\in \mathcal{J}$ define
    \begin{equation}\label{xi}
      \xi_m := \sqrt{m!}\prod_{k=1}^\infty H_{ m_{k}}\big(\chi_{k} \big),
    \end{equation}
    where $H_n$ is the $n$-{\rm th} Hermite polynomial, $m!:=\prod_{k=1}^\infty m_{k}!  \;$ and 
    \begin{equation}\label{chi_def}
    \chi_{k}:= \int_0^T e_k(s) dw_k(s) \, \sim N(0,1).
    \end{equation}
    It is well-known that the family $\Big\{ \xi_m : m\in  \mathcal{J} \Big\} $ is a complete orthonormal basis of $L^2(\Omega, \mathcal{F},\Pb)$ (cf. 
    \cite[Theorem 4.2]{lototsky2006stochastic}). Indeed, for $n\in\IN$ denoted by 
    $\mathcal{H}_n$ as the closed subspace of  $L^2(\Omega, \mathcal{F},\Pb)$
    generated by the family of random variables $\xi_m$ with $|m|=n$, then 
    \begin{equation}\label{cameron-martin-decomposition}
        L^2(\Omega, \mathcal{F},\Pb) = \bigoplus_{n=0}^\infty \mathcal{H}_n.
    \end{equation}
    We will also refer to $ \mathcal{H}_n$ as the $n$-th order Wiener chaos. Then, any random variable $u\in L^2(\Omega) $ can be written as
    \begin{equation}\label{cameron-martin}
        u= \sum_{m\in \mathcal{J}}  u_m \xi_m,
    \end{equation}
    where $\xi_m$ is given by \eqref{xi} and $u_m = \E(u\, \xi_m)$.\\
    We now introduce the notion of Malliavin derivative and its adjoint operator: 
    the Skorohod integral in our particular setting.\\
    A random variable $F\in L^2(\Omega)$ is smooth if there exists
    $h_1,h_2,\dots, h_n\in \mathbb{H}$ such that 
    \begin{equation}\label{smooth-RV}
        F=\psi\big(W(h_1),\ldots,W(h_n)\big),
    \end{equation}
    where the function $\psi:\R^n \rightarrow \R$ is infinitely differentiable with 
    all derivatives polynomial bounded. We denote the set of all smooth random 
    variables given by \eqref{smooth-RV} as $\mathcal{S}$.\\
    We have the following definition.
    \begin{definition}
    Let $F$ be the smooth random variable defined by \eqref{smooth-RV}. The
    {\rm Malliavin derivative} of $F$ is given by
        \begin{equation*}
            DF= \sum_{i=1}^n \partial_i \psi\big(W(h_1),\ldots,W(h_n)\big) 
            \, h_i.
        \end{equation*}
    \end{definition}
    Observe, that the Malliavin derivative $DF$ is an $\mathbb{H}$-valued random
    variable. In our setting, we write $D_tF$ for the realization of the function 
    $DF$ at the point $t\in [0,T]$:
    \begin{equation*}
        D_tF= \sum_{i=1}^n \partial_i  \psi\big(W(h_1),\ldots,W(h_n)\big) \, h_i(t).
    \end{equation*}
    The Malliavin derivative is a closable operator  from $L^p(\Omega)$ to
    $L^p(\Omega;\mathbb{H})$, for all $p\ge 1$; thus we consider its closure as 
    our derivative operator, which we will still denote
    by $D$.
    
    The space $\mathbb{D}^{1,p}$ denotes the completion of $\mathcal{S}$ with 
    respect to the norm %
    \begin{equation*}
        \|F \|_{1,p}^p  := \E\big( |F|^p \big) +  
        \E\big[\|D F \|_{\mathbb{H}}^p \big].
    \end{equation*}
    In our setting, we consider the case $p=2$.\\
    The Malliavin derivative, at point $s\in [0,T]$, of the random variable 
$\xi_m$ given by \eqref{xi} is
        \begin{equation}\label{xi-derivative}
       D_s\xi_m = \sum_{k=1}^\infty \sqrt{m_{k}}\, e_k(s)\, 
            \xi_{m^-(k)}\, y_l,
        \end{equation}
    where
        \begin{equation}\label{m-i}
          \big(m^-(k)\big)_{i}:= {\begin{cases} \max\big\{m_{k}-1\, ,\, 0\big\}
            \quad & \mbox{ if } i=k,\\ m_{k} & {\mbox{otherwise.}}
            \end{cases}}
        \end{equation}
    We now review the general definition of the adjoint of operator $D$. Since $D$ is 
    a closed and unbounded operator from its domain, which is the dense-closed subspace 
    $\mathbb{D}^{1,2}$ of $L^2(\Omega)$, to $L^2(\Omega;\mathbb{H})$. Then, we can define 
    an adjoint operator, denoted by $\delta$, which maps a dense subspace of 
    $L^2(\Omega;\mathbb{H})$ to $L^2(\Omega)$. More formally, denote by $\delta$ the adjoint of $D$. The domain of $\delta$, 
    denoted by $Dom (\delta)$ is the set of $\mathbb{H}$-valued square-integrable
    random variables $u\in  L^2(\Omega;\mathbb{H})$ such that 
     \begin{equation*}
         \Big |\E\big(\langle DF, u \rangle_{\mathbb{H}}  \big) \Big|\le C
         \left(\E | F |^{2} \right)^{1/2},
     \end{equation*}
    for all $F\in \mathbb{D}^{1,2}$  and $C$ is a constant depending on $u$.\\ 
    Moreover, if $u\in Dom (\delta)$, then $\delta(u)$ is the element of 
    $L^2(\Omega)$ characterized by
    \begin{equation}\label{def-delta}
        \E\big(F\delta(u) \big)= \E\big(\langle DF, u \rangle_{\mathbb{H}}\big),
    \end{equation}
    for any $F \in \mathbb{D}^{1,2}$.
    \begin{remark}\label{remark-ito-skorohod}
    The operator $\delta$ is called the divergence operator or the Skorohod integral.
    In \cite{nualart1988stochastic} has been proved that, under some not very restricted
    assumptions,  the It\^o stochastic integral is a particular case of the Skorohod
    integral. Indeed, the only additional hypothesis is that $u$ is adapted to the 
    filtration generated by the noise $W_t$. If the framework of our work does not 
    satisfy such hypotheses, then we cannot write the It\^o stochastic integral as 
    the divergence operator (see Section 1.3 in \cite{nualart2006malliavin}). 
    As an application of a result given by Lototsky 
    and Rozovzky  in \cite{lototsky2006stochastic} (see Lemma 5.1 therein) we have
        \begin{equation}\label{diver_Der}
            \E\big(U(t) \xi_m\big) = \E \int_0^t
            u(s) (D_s \xi_m) ds,
        \end{equation}
    where 
    $$
    U(t)= \int_0^t u(s) dB_s
    $$
  We will use 
    this  formula \eqref{diver_Der} to obtain the 
    WCE of the  It\^o SDEs.
    \end{remark}
    \section{Wiener Chaos expansion for the diffusion bridge} \label{WienerChaos-BM}
    In this section, we present the propagator system, that is, a system of ordinary 
    differential equations (ODEs) whose solutions are the coefficients of the WCE 
    \eqref{cameron-martin}. For further reading on the propagator for general SDE 
    and stochastic partial differential equations, we refer to \cite{lototsky2017stochastic}.
    
    First, we present a procedure to obtain the propagator for the 
    SDE \eqref{SDE1}. We know that the solution $\bm X$ belongs to the space 
    $L^2(\Omega; [0,T])$ then, we can write the solution as in \eqref{cameron-martin}:
    \begin{equation}\label{X-WCE}
      X(t)=\sum_{m\in\mathcal{J}} X_m(t) \xi_m, 
    \end{equation}
    where $ X_m(t)=\E\big( X(t) \xi_m \big) $. Thus, using the integral representation
    of the SDE \eqref{integral_SDE}
        \begin{align}\label{xm1}
            X_m(t) 
                &=  \E\big( X(t) \xi_m \big) \nonumber\\
                &= \E\Bigg[\Big(x_0 + \int_0^t  f(X(s))ds + \int_0^t   \sigma(X(s))dB(s)  \Big)  
                \xi_m \Bigg] \nonumber\\
                &= x_0 \I_{\{|m|= 0\} }  + \int_0^t  \E\big[f\big(X(s)\big) \xi_m\big] ds 
                + \E\Bigg[\Bigg(\int_0^t   \sigma(X(s))dB(s) \Bigg)   \xi_m \Bigg]\nonumber\\
                &= x_0 \I_{\{|m|= 0\} }  + \int_0^t  f_m(X(s)) ds 
                + \E\Bigg[\Bigg(\int_0^t   \sigma(X(s))dB(s) \Bigg)   \xi_m \Bigg].
        \end{align}
    Here $f_m(X):=\E\big[f\big(X\big) \xi_m\big] $. For the last integral in \eqref{xm1}, we will use the Remark \ref{remark-ito-skorohod}. Indeed, we can write the following
    \begin{equation*}
         \int_0^t \sigma\big(X(s)\big)dB(s) =
         \delta\Big( \sigma\big(X(\cdot)\big) \I_{[0,t]}(\cdot)\Big).
    \end{equation*}
    Thus, using this expression, \eqref{xi-derivative} and \eqref{diver_Der} we get
    \begin{align}\label{A-term}
        \E\bigg[ \Big( \int_0^t \sigma\big(X(s)\big)dB(s) \Big) \xi_m \bigg] 
        &= 
            \E\Bigg( \delta\Big( \sigma\big(X(\cdot)\big) \I_{[0,t]}(\cdot)\Big) 
            \xi_m \Bigg)\nonumber\\
        &= 
            \E\Bigg(\Big\langle \sigma\big(X(\cdot)\big) \I_{[0,t]}(\cdot) ,
            D \xi_m \Big\rangle_{\mathbb{H}} \Big)\nonumber\\
        &= 
            \E\bigg[  \int_0^t  \sigma\big(X(s)\big)\,D \xi_m\,  ds  \bigg] \nonumber\\
        &= \sum_{i=1}^\infty \sqrt{m_{i}} \int_0^t e_i(s) \E\Big[ \xi_{m^-(i)}
       \sigma\big(X(s)\big)\Big] ds \nonumber\\
        &= \sum_{i=1}^\infty \sqrt{m_{i}}  \int_0^t e_i(s) \sigma_{m^-(i)}
        \big(X(s)\big) ds.
    \end{align}
    Plugging together \eqref{xm1} - \eqref{A-term} we have
    \small{
        \begin{align}\label{xm-equation}
            X_m(t) &=  x_0  \I_{\{|m|= 0\} } + \int_0^t  f_m(X(s)) ds 
            + \sum_{i=1}^\infty \sqrt{m_{i}} \int_0^t e_i(s) 
            \sigma_{m^-(i)}\big(X(s)\big) ds.
        \end{align}}
    
    \begin{remark}
        At this point,  we could try to include the final condition $X_T=x_T$ into 
    the propagator system \eqref{xm-equation}; that is, using \eqref{X-WCE} 
    we obtain
    \begin{align}
        x_T &=X(T)= \sum_{m \in \mathcal{J}} X_m(T)\xi_m= x_T \I_{|m|=0} 
                + \sum_{\{m \in \mathcal{J}: |m|\ge 1\}} X_m(T)  \xi_m.
    \end{align}
    From this, we deduce that the condition for the propagator $X_m$ at the point 
    $t=T$ would be
    $$
    X_m(T)=  x_T \I_{|m|=0}.
    $$
    This implies that the propagator for $|m|=0$ is a  first-order differential 
    equations with boundary conditions, which have no solution. Therefore, it is 
    necessary to use another approach to simulate diffusion bridges with the WCE.
    \end{remark}
    
\subsection{A representation for the diffusion bridge}
    In this subsection, we use the particular representation \eqref{Yt-eq} for a 
    diffusion bridge, this representation allows us to avoid the problem of solving 
    a first-order differential equation with boundary conditions. In addition, 
    we deduce a WCE for this representation which permits us to simulate a bridge for 
    the SDE \eqref{SDE1}.
    
    Consider the  process $\bm Y$ given by \eqref{Yt-eq}:

 \begin{equation*}
    dY(t) =  (\theta-Y(t)) \frac{1}{T-t} dt + dX(t),\qquad Y(0)=\eta,
    \end{equation*}

with solution given by
 \begin{equation}\label{sol_bri}
    Y(t) = \eta+(\theta-\eta) \frac{t}{T}+(T-t) 
                \int_{0}^{t} \frac{d X(s)}{T-s}.
    \end{equation}

\begin{remark}
Note that if we have observations at high frequency of the path $X(\cdot)$, thus  we could approximate the integral, in \eqref{sol_bri}, numerically; however, this is not the case.
\end{remark}
  Observe that the process $\bm Y$ satisfies $Y(0)=\eta$,\, and $Y(T)=\theta$.  Thus, the process $\bm Y$ is a proposal bridge process from $(0, \eta$) to 
    $(T, \theta)$.
    
    Since the SDE $\bm X$ belongs to $L^2(\Omega; \mathbb{H})$, $\bm Y\in L^2(\Omega; 
    \mathbb{H})$ too. Then, we write the WCE of $\bm Y$:
    \begin{equation}\label{Yt-propagator}
        Y(t)=\sum_{m \in \mathcal{J}} Y_m(t) \xi_m,
        \quad \mbox{ with }\quad Y_m(t) = \E(Y(t) \xi_m).
    \end{equation}
    The propagator $Y_m(t)$ solves 
    \begin{equation*}
        \begin{aligned}
            Y_{m}(t) 
                &=\mathbb{E}\left[\left(\eta+(\theta-\eta) \frac{t}{T}+(T-t) 
                \int_{0}^{t} \frac{d X(s)}{T-s}\right) \xi_{m}\right] \\
                &=\left(\eta-(\theta-\eta) \frac{t}{T}\right) 
                    I_{\{|m|=0\}}+(T-t) \mathbb{E}\left[\left(\int_{0}^{t} 
                    \frac{d X(s)}{T-s}\right) \xi_{m}\right].
        \end{aligned}
    \end{equation*}
    For the term given by the stochastic integral with respect to $\bm X$. We take 
    a partition of the interval $[0, t]$ in $n$ uniform pieces, then we can write
    \begin{equation*}
        \begin{aligned}
            \mathbb{E}\left[\left(\int_{0}^{t} \frac{d X(s)}{T-s}\right) 
            \xi_{m}\right] &=\mathbb{E}\left[\lim_{n \rightarrow \infty} 
            \sum_{i=0}^{n} \frac{1}{T-t_{i}}
            \left(X\left(t_{i}\right)-X\left(t_{i-1}\right)\right) \xi_{m}\right]
            \\
            &=\lim_{n \rightarrow \infty} \sum_{i=0}^{n} \frac{1}{T-t_{i}} 
            \mathbb{E}\left[\left(X\left(t_{i}\right)-X\left(t_{i-1}\right)\right) 
            \xi_{m}\right] \\
            &=\lim_{n \rightarrow \infty} \sum_{i=0}^{n} 
            \frac{1}{T-t_{i}}\left[\mathbb{E}\left(X\left(t_{i}\right) \xi_{m}\right)
            -\mathbb{E}\left(X\left(t_{i-1}\right) \xi_{m}\right)\right] \\
            &=\lim_{n \rightarrow \infty} \sum_{i=0}^{n} 
            \frac{1}{T-t_{i}}\left[X_{m}\left(t_{i}\right)-X_{m}\left(t_{i-1}\right)
            \right]\\
            &= \int_{0}^{t} \frac{d X_m(s)}{T-s},
        \end{aligned}
    \end{equation*}
    thus
    \begin{align}\label{Ym}
        Y_m(t) & =  \left(\eta+(\theta-\eta) \frac{t}{T}\right) I_{\{|m|=0\}}+(T-t) 
        \int_{0}^{t} \frac{d X_m(s)}{T-s},
    \end{align}
    where $X_m(t)=\E(X_t \xi_m)$ is the propagator for $X$ given $m\in\mathcal{J}$. 
    
    We use equations \eqref{xm-equation}, \eqref{Yt-propagator} and \eqref{Ym} to 
    simulate a bridge for the diffusion $\bm X$. 
    
We rewrite Equation (17) as an ODE system,
\begin{equation*}
    \begin{aligned}
        \frac{dX_m}{dt} &= f_m(X(t)) + \sum^{\infty}_{i=1}\sqrt{m_i} e_i(t)\sigma_{m^{-}(i)}(X(t))\\
        X_{m}(0)&= x_0 I_{\{|m|=0\}}.
    \end{aligned}
\end{equation*}
We use the solution $X_m(t)$ to obtain $Y_m(t)$ in \eqref{Ym}. Note that, for all 
$t \in (0,T]$ and $s \in(0,t)$, we have $\frac{T-t}{T-s}< 1$. Consequently, the value of the stochastic integral presented in \eqref{Ym} does not diverge.

    \subsection{An approximation to the  Wiener-Chaos Expansion} 
    \label{sect_approx_WCE} In this subsection, we truncate the WCE to construct an approximation for the diffusion and therefore for the 
    diffusion bridge.
    
    Define 
        \begin{equation}\label{trunc_J}
            \mathcal{J}_{p,L} := \Big\{ m = (m_{k})_{1\le k\le L},\,  
            \,|m| \le p \Big\}.
        \end{equation}
    We denote
    \begin{equation}\label{trunc_WCE_X}
        X_t^{p, L}=\sum_{m \in \mathcal{J}_{p,L}} X_m (t)\, \xi_m,
    \end{equation}
    with coefficient functions $X_m(t) = \E(X(t) \, \xi_m)$.
    
    Thus, we have
    \begin{equation*}
        X_t^{p, L}\quad \longrightarrow X(t) =\sum_{m \in \mathcal{J}} X_m (t) 
        \quad \mbox{ in } L^2, \mbox{ when } p, L\rightarrow \infty \, .
    \end{equation*}
    The mean-square error of the truncated WCE at the terminal time 
    $T$ (that is $X_{T}^{p, L}$), depends on the order $p$ and the maximum length $L$ 
    of the expansion (i.e., the number $L$ of included basis functions $e_i(\cdot)$) 
    and can be estimated by (see \cite{huschto2019asymptotic} for instance)
        \begin{equation}\label{error-WCE}
            \mathbb{E}\left[\left|X(T)-X_{T}^{p, L}\right|^2\right] 
            \le C_1\left(1+x_0^2\right) e^{\left(C_1+\kappa^2\right) T} 
            \frac{\left(\kappa^2 T\right)^{p+1}}{(p+1)!}+C_2\left(1+x_0^2
            \right)    \frac{T^4}{L} e^{C_2 T},
        \end{equation}
    with constants $C_1, C_2$ depending on $K$ and $T$, and $\kappa$ depending only 
    on $K$. Here $K$ is the Lipchitz constant for the coefficients of the SDE.
    
    Thus, we need to take $p$ and $L$ large enough such that the error is small 
    enough to obtain a good approximation for the diffusion bridge.
    \section{Validation of the method} \label{validation_WCE}
        Let us begin this section, with the presentation of a generic algorithm for simulating diffusion bridges using our WCE method. We illustrate this algorithm with the Ornstein-Uhlenbeck (OU) process. Furthermore, to validate our method, we compare the numerical results offered using the WCE with the exact method for the OU process (cf. \cite{bladtSoren},\cite{bladtSorenc}).
\\\\
    The generic Algorithm \ref{al-WCE} is as follows. 
    
    \begin{algorithm}[H]
        \begin{algorithmic}[1]
            \STATE Choose the algorithm's parameters:  
                \begin{itemize}
                    \item A basis $\{e_m(x)\}_{m\ge 1}$ for $L^2(0,T)$ and its truncation
                    $\{e_m(x)\}_{m\ge 1}^L$.
                    \item A time step $\Delta$ such that $K \Delta =T$.
                    \item Choose an initial value for $X$, namely, $x_0$.
                    \item Define the truncation of $\mathcal{J}$, meaning $\mathcal{J}_{p,L}$ (see \eqref{trunc_J}),
                        where $p$ is the highest order of polynomials to be considered and $L$ is the length of the vector.
                    \item Fix a finite number of elements of $\mathcal{J}_{p,L}$, i,e,
                        $\{m^i \}_{i=1}^L$ with each $m^i\in \mathcal{J}_{p,L}$.
                \end{itemize}
            \STATE For each $i=1,\ldots,pL$, solve the propagator \eqref{xm-equation}
            for the $m^i\in \mathcal{J}_{p,L}$ already fixed on the time interval $[0,T]$ 
            with its corresponding initial condition $X^i(0)=x_0^i$.
            \STATE Generate
            the Gaussian random variables $\chi_{m^i}$,  $m^i\in \mathcal{J}_{p,L}$, for $i=1,\ldots,pL$.
            \STATE Obtain the approximate solution of the SDE \eqref{SDE1} \newline
            \STATE Construct the  approximation of the diffusion bridge using \eqref{Yt-propagator} and \eqref{Ym}.
        \end{algorithmic}
    \caption{Simulation of diffusion bridges using the WCE.}\label{al-WCE}
    \end{algorithm}
\subsection{OU bridge with the WCE}
    Here, we use Algorithm \ref{al-WCE} to illustrate the WCE method to simulate 
    diffusion bridges of an OU process. Set the Hilbert space $\mathcal{H}=L^2(0,T)$ 
    and define a basis in this space 
    as
        \begin{equation*}
            e_j(t):=\sqrt{\frac{2}{T}} \sin\big(j\pi\, \tfrac{t}{T}\big),
                \qquad t\in [0,T].
        \end{equation*}
    Consider the diffusion process $\bm X$ solution of the SDE 
        \begin{align}\label{OU-SDE}
            dX(t) &= -a X(t) dt + \sigma dB(t),\\    
            X(0)&= x_0,  \nonumber
        \end{align}
    where $a,\sigma>0$, with its integral representation 
        \begin{align}\label{OU-integral}
            X(t) &= x_0 -a \int_0^t X(s) ds + \sigma \int_0^t  dB(s). 
        \end{align}
    With suitable assumptions we have that $\bm X\in L^2(\Omega, \mathcal{F},\Pb)$, 
    this implies that we can write
        \begin{equation*}
            X(t)= \sum_{m \in \mathcal{J}} X_m(t) \xi_m, 
                \quad \mbox{with }  
            X_m(t)= \E(X(t)\, \xi_m),
        \end{equation*}
    thus, using the one-dimensional version of \eqref{diver_Der} we get
    \begin{align}\label{propagator_OU}
        X_m(t)
            &= \E\Big[\Big(x_0  -a \int_0^t X(s) ds 
                + \sigma \int_0^t  dB(s) \Big)\xi_m\Big] \nonumber\\   
            &=  x_0\I_{|m|=0}   -a \int_0^t  \E\Big(X(s)\xi_m \Big)  ds 
                + \sigma \E\Big[ \Big(\int_0^t  dB(s)\Big) \xi_m\Big] \nonumber\\
            &= x_0\I_{|m|=0}  - a \int_0^t  X_m(s) ds 
                + \sigma \E\Big[ \big\langle \I_{[0,t]}(\cdot), D\xi_{m}
                \big\rangle_{\mathcal{H}} \Big] \nonumber\\
            &= x_0\I_{|m|=0}  - a \int_0^t  X_m(s) ds 
                + \sigma \int_0^t \sum_{i=1}^\infty \sqrt{m_{i}} \, e_i(s) \, 
                    \E \big( \xi_{m^-(i)} \big) ds\nonumber\\
            &=x_0\I_{|m|=0}  - a \int_0^t  X_m(s) ds 
                + \sigma\sum_{i=1}^\infty \int_0^t e_i(s) ds \I_{\{i=1\} } \,\, 
                \I_{\{|m^-(i)|=0\}}.
    \end{align}
    From this expression, we now calculate the propagator for the following three 
    cases $|m|=0$, $|m|=1$, and $|m|>1$.
    
    For $|\alpha|=0$ denote the propagator as $ X_0:= X_{|m|=0}$, observe that 
    in this case the indicator function $\I_{|m^-(i,j)|=0}$ is equal to zero, 
    then from \eqref{propagator_OU} we have the ordinary differential equation
    \begin{align*}
    \frac{dX_0(t)}{dt}
    &= - a   X_0(t),
    \end{align*}
    with initial condition $x_0$. This means that the solution $X_0(t)$ is
    \begin{align*}
     X_0(t)
    &= x_0 \exp(- at).
    \end{align*}
    For the case, $|m|=1$ denote the propagator as $ X_1= X_{|m|=1}$. 
    Then, from \eqref{propagator_OU} and given that $|m|=1$ then we have  
        $\sqrt{m_{i}}=1 $ for some $i$, thus we can write for such $i$
    \begin{align}
        \frac{dX_1(t)}{dt}
            &= - a X_1(t)  + \sigma  \sum_{i=1}^\infty \sqrt{m_{i}}\,\,
            e_i(t)  \I_{\{i=1\} } \,\, \I_{\{|m^-(i)|=0\}}\nonumber\\
            &= - a   X_1(t) + \sigma e_i(t), \nonumber
    \end{align}
    where the ODE has zero initial condition. Thus, we have that 
    \begin{align}\label{m1_OU}
     X_{1_i}(t) =  \sigma  e^{-at} \int_0^t e^{as} e_i(s) ds.
    \end{align}
    For the case, $|m|>1$ we have the ODE with zero initial condition 
    \begin{align*}
    \frac{dX_m(t)}{dt}
    &= - a   X_m(t), 
    \end{align*}
    with the trivial solution $X_m(t) \equiv 0$. This implies that the WCE for the OU process could be written as
    \begin{align*}
    X(t) = \sum_{m \in \mathcal{J}} X_m(t) \xi_m = X_0(t) 
    + \sum_{k =1}^\infty  
    X_{1_k}(t) \xi_{m_k},
    \end{align*}
    with $X_{1_k}$ given by \eqref{m1_OU} and
    \begin{equation*}
        \xi_{m_k} :=  \int_0^T e_k(s) dB(s).
    \end{equation*}

    Recall that the solution of the OU is
    \begin{equation*}
        X(t)= x_0 e^{-at} + \sigma e^{-at} \int_0^t  e^{as} dB(s),
    \end{equation*}
    and since the Wiener integral is Gaussian with zero mean and variance given by 
    the It\^o isometry, we could rewrite the solution as
        \begin{align*}
            X(t)= x^0 e^{-at} + \frac{\sigma}{\sqrt{2a}}   B(1-e^{-2at}).   
        \end{align*}
    In addition, we know that the Brownian motion has a WCE given by
        \begin{equation*}
            B(t)= \sum_{k\ge 1} \Big(  \int_0^t e_k(s) ds \Big)\, \xi_{1_k},
        \end{equation*}
    where 
    \begin{equation*}
        \xi_{1_k}:=\int_0^T e_k(s) dB(s).
    \end{equation*}
    Using the previous, we define our approximation as 
    \begin{align*}
    X^L(t)  = X_0(t) +   \hat X_1(t)\diamond \hat\xi_1,
    \end{align*}
    where we have defined
        \begin{equation*}
            \hat X_1(t)\diamond \hat\xi_1 := \sigma e^{-at} \sum_{k=1}^L \xi_{1_k} 
            \int_0^t e^{as} e_k(s) ds.
        \end{equation*}
        
    \subsection{Simulation study}
    Now, we present a simulation study to calibrate the approximation 
    method of diffusion bridges based on the WCE presented in this work. 
    We compare the diffusion bridges obtained by our approximate method with the 
    following exact method (see \cite{bladtSoren},\cite{bladtSorenc}) to construct the OU
    bridge.
    
    Let $0=t_0<t_1,\ldots<t_n=T$ be a discretization of time interval $[0,T]$ 
    and $Y(0)=\eta$. If we generate 
        \begin{equation*}
            Y(t_i)=e^{-a(t_i-t_{i-1})}Y(t_{i-1})+W_i,
        \end{equation*}
    where $W_i\sim N\left(0,\frac{\sigma^2(1-e^{2a(t_i-t_{i-1})})}{2a}\right)$, 
    for $i=1,\ldots,n.$ Then,
        \begin{equation*}
            Z(t_i)=Y(t_i)+(\theta-X(T))\big(\frac{e^{at_i}-e^{-at_i}}
            {e^{aT}-e^{-aT}}\big),
            \hspace{1cm}i=0,1,\ldots,n,
        \end{equation*}
    is an OU bridge from $(0,\eta)$ to $(T,\theta)$.
    
    To calibrate the Wiener chaos approximation method, we consider four examples of OU bridges from $(0,0)$ to $(1,0),(1,1)$, and $(1,2)$, and from $(0,0.8)$ to 
    $(1,0.5)$. For each bridge, we simulated $1,000$ paths and we conducted tests with several combinations of parameters, but to illustrate our method, we only present the results with parameters $a = 0.5$ and $\sigma =1.0$.
    
\begin{figure}[H]
        \centering
        \includegraphics[width=1.0\textwidth]{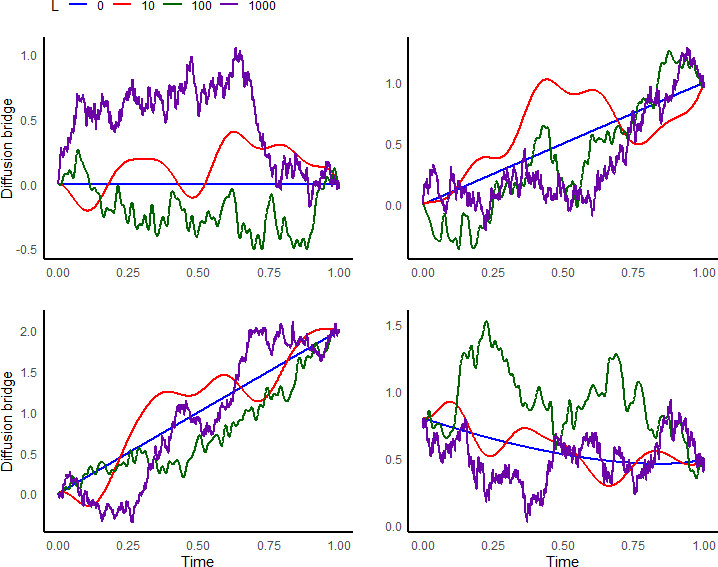}
        \caption{%
       OU bridges with different values of $L$ from $a)$ $(0,0)$ to $(1,0)$, 
        $b)$ $(0,0)$ to $(1,1)$, $c)$ $(0,0)$ to $(1,2)$, and 
        $d)$ $(0,0.8)$ to $(1,0.5)$  with parameters 
        $a = 0.5$, $\theta = 1.0$.}
        \label{fig:OU_bridges1}
    \end{figure}

    In Figure \ref{fig:OU_bridges1}, we present OU bridges 
    for the four cases. The bridges are simulated using the Wiener chaos approximation method. Figure \ref{fig:OU_bridges1}  illustrates, on one hand, the deterministic part of the process, and on the other hand, the effect of increasing the number $L$ (the level of truncation) is included in the approximation of the WCE. We ran the simulations considering $p=12$ (see \eqref{error-WCE}) and $L=0,10,100, 1000$, where $L$ is the level of truncation.
    
\begin{figure}[H]
     \centering
      \subfloat[]{
       \label{fig:ou1}
        \includegraphics[width=0.48\textwidth,keepaspectratio]{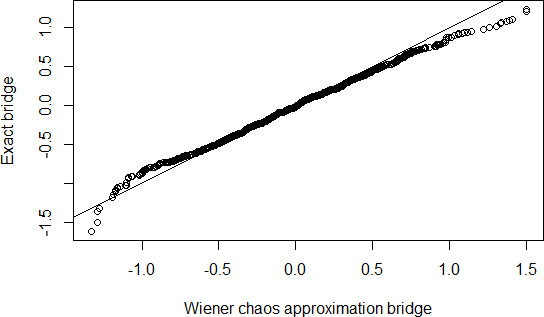}}
      \subfloat[]{
       \label{fig:ou2}
        \includegraphics[width=0.48\textwidth,keepaspectratio]{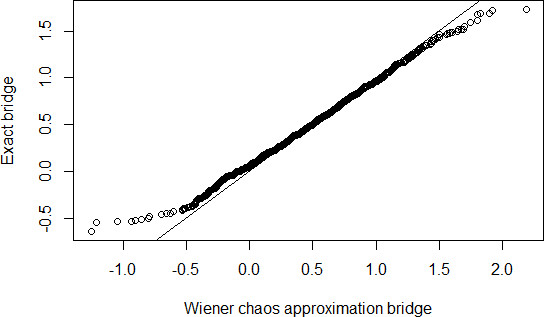}}
      \hfill
      \subfloat[]{
       \label{fig:ou3}
        \includegraphics[width=0.48\textwidth,keepaspectratio]{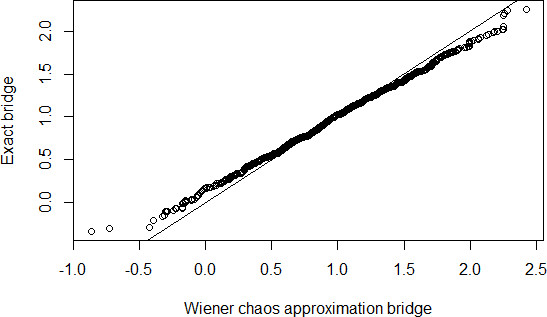}}
         \subfloat[]{
       \label{fig:ou4}
        \includegraphics[width=0.48\textwidth,keepaspectratio]{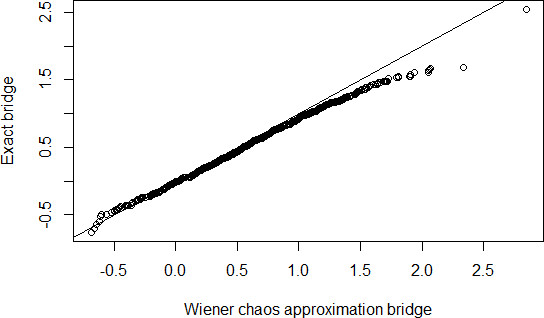}}
     \caption{Comparison of  the empirical distributions between WCE and exact method 
     at time 0.5 based on 1,000 realizations of Ornstein-Uhlenbeck bridges from  
     $(a)$ $(0,0)$ to $(1,0)$, $(b)$ from $(0,0)$ to $(1,1)$, $(c)$ from $(0,0)$ to 
     $(1,2)$ and $(d)$ from $(0,0.8)$ to $(1,0.5)$.}
     \label{f:q-q}
    \end{figure}

   Figure \ref{f:q-q} presents a QQ-plot, which is a scatter plot that illustrates the relationship between two sets of quantiles. In this figure, we compare the quantiles derived from our approximate method with those from the empirical distribution obtained via the exact algorithm. In all four cases, the points lie close to a straight line, suggesting that both sets of quantiles probably come from the same distribution.
    
    Moreover, we apply the Kolmogorov-Smirnov (KS) test, a nonparametric test that compares the cumulative distributions of two data sets, to compare samples from our method with those from the exact method. A small p-value from the KS test (p-value < 0.05) indicates that we can reject the null hypothesis, leading us to conclude that the two sets of points are likely drawn from distinct distributions. The results from the KS test, as displayed in Table \ref{tb:CPU-values}, consistently exceed the conventional significance threshold of 0.05. This provides strong evidence that our method and the exact method share the same empirical distribution. This indicates that our proposed method for constructing diffusion bridges demonstrates strong performance.
        \begin{table}[H]
    \centering
    \begin{tabular}{lcrrrr}
      \hline
     $(0,\eta) \to (T, \theta)$ & L & CPU  &KS test& \\ 
     & (truncation) & (sec.) &  (p-value) \\ 
      \hline
      $(0,0) \to (1,0)$ & 100 & 0.0028 & 0.37 \\ 
      &  1000 & 0.0273 & 0.1995 \\
        &  10000 & 0.2239 & 0.2877 \\
      $(0,0) \to (1,1)$ & 100 & 0.0026 & 0.9688 \\ 
      &  1000 & 0.0219 & 0.37 \\
        &  10000 & 0.3188 & 0.2193\\
      $(0,0) \to (1,2)$ & 100 & 0.0036 & 0.0546 \\ 
      &  1000 & 0.0292 & 0.2877 \\
        &  10000 & 0.2842 & 0.2406 \\
      $(0,0.8) \to (1,0.5)$ & 100 & 0.0026 & 0.1082 \\ 
      &  1000 & 0.0242 & 0.0971\\
        & 10000 & 0.1056 & 0.0954 \\
     \hline
    \end{tabular}
        \caption{%
            In this table, we are presenting several information. First, we report the CPU execution time average for construction one diffusion bridge using the WCE with different levels of truncations (recall the term $L$ in \eqref{error-WCE}). We also present the results of Kolmogorov-Smirnov test for the WCE method and the exact method (EM). We simulate 1000 Ornstein-Uhlenbeck 
            bridges. We vary the initial and endpoints.
        }
        \label{tb:CPU-values}
    \end{table}

    The simulations were performed on an ASUS H110M-A/M.2 (Intel Core i5 7400) quad-core 3.00
    GHz with a memory RAM of 8 GB.  The CPU execution time applied to the Ornstein-Uhlenbeck
    bridges using the Wiener approximation method is reported in Table \ref{tb:CPU-values}.
    These simulations were made in Julia and R language. We can observe that when we use 100 
    or 1000 Brownian motions, the CPU execution time is slightly larger in comparison with 
    the approximate method from  \cite{bladtSoren}, and \cite{bladtSorenc}. Indeed, in these cases, the CPU execution
    time is around 0.5-2.0 CPU seconds faster. In the case of 10000 Brownian motions, the CPU
    time of the Wiener approximation method was around 285 seconds, giving an accurate form 
    to construct a diffusion bridge.

    As the WCE is an approximation method, one goal is to accurately determine the minimum number of Brownian motions needed to achieve a high-quality approximation (see the term $L$ in equation \eqref{error-WCE}).
     To determine this number, we performed the KS test, using as a reference the exact method, to determine the point at which $L$ is sufficiently large to yield a p-value greater than $0.05$. The results are reported in Table \ref{tb:min-CPU-values}. 
    From this results, we can conclude that, for the OU process, the WCE method with a minimum level of $L$ can build diffusion bridges. 
    \begin{table}[H]
    \centering
    \begin{tabular}{lcrr}
      \hline
     $(0,\eta) \to (T, \theta)$ &  Min. L & CPU & KS test \\ 
     & (truncation) & (sec.) &(p-value)\\ 
      \hline
      $(0,0) \to (1,0)$ & 5 & 0.001460046&  0.1640792\\ 
      $(0,0) \to (1,1)$ & 10 & 0.001340425 & 0.1338343\\ 
      $(0,0) \to (1,2)$ & 5 & 0.001314852 & 0.06919033\\ 
      $(0,0.8) \to (1,0.5)$ & 25 & 0.001885219 & 0.07762146\\ 
     \hline
    \end{tabular}
        \caption{%
            Minimum value of $L$  used to approximate efficiently Ornstein-Uhlenbeck bridges  by using the             Wiener chaos expansion method. The results of the Kolmogorov-Smirnov test were obtained by implementing 1000 simulations for both the exact method and the WCE method. We also report the CPU execution time average for construction one diffusion bridge with the reported minimum BM.
        }
        \label{tb:min-CPU-values}
    \end{table}

    \section{Numerical Validation}\label{numerics}
  In this section, we compare our method with existing numerical methods for approximating diffusion bridges for the OU process and Geometric Brownian motion. Additionally, we apply Algorithm \ref{al-WCE} to two other numerical examples: the affine logistic SDE and the stochastic protein kinetic differential equation. For all these examples, we solve numerically  their propagator. We used the package \textit{DifferentialEquations.jl} from 
    \textit{Julia}. This package provides recent solvers for unknown stiffness problems. 
    We applied the \textit{AutoTsit5(Rosenbrock23())} method, getting a stable, and 
    efficient simulation for large systems ($> 1000$ ODEs).\\
    
    For the simulations in this section, we consider $T=1$ and fixed the Hilbert space $L^2(0,1)$ and take 
    the basis defined by
    \begin{equation*}
        e_j(t)=\sqrt{\frac{2}{T}} \sin\big(j\pi\, \tfrac{t}{T}\big),\qquad t\in [0,T].
    \end{equation*}

 \subsection{Ornstein Uhlenbeck process}

    In this subsection, we thoroughly examine the OU bridge once again. However, our primary focus is now centered upon the comparison between the outcomes of the WCE method and two other alternative methods. Indeed, we have simulated the OU bridges with the methods called Doob's $h$-transform (Doob's $h$.) and  Bladt and S\o rensen's (B\&S) described in \ref{other-method-app}. 
    
    Based on the findings presented in Table \ref{tb:Comparison_OU}, it is evident that the WCE method aligns perfectly with Doob's $h$-transform approach. This conclusion is drawn from the p-values greater than $0.05$ obtained through the Kolmogorov-Smirnov test, which remain consistent across all the simulated scenarios and cases analyzed.

    \begin{table}[H]
    \centering
    \begin{tabular}{lcrrrr}
      \hline
     $(0,\eta) \to (T, \theta)$ & L  & Doob's $h$ &  B\&S \\ 
     & (truncation) &  (p-value) & (p-value)\\ 
      \hline
      $(0,0) \to (1,0)$ & 100 &  0.3699 & 0.2633\\ 
      &  1000 & 0.341 & 0.0546\\
        &  10000 &   0.2406 & 0.9134\\
      $(0,0) \to (1,1)$ & 100 &  0.6099 & 0.0691\\ 
      &  1000 & 0.7943 & 0.0869\\
        &  10000 &  0.0869 & 0.0546\\
      $(0,0) \to (1,2)$ & 100 &  0.5726 & 0.0615\\ 
      &  1000 & 0.1483 & 0.0530\\
        &  10000& 0.1811 & 0.0004\\
      $(0,0.8) \to (1,0.5)$ & 100 &   0.2917 & 0.2574\\ 
      &  1000 &  0.2262 & 0.0693\\
        & 10000 &  0.0546 & 0.1205\\
     \hline
    \end{tabular}
        \caption{%
            In this table, we are presenting the results of the Kolmogorov-Smirnov test for the WCE method and the two  methods: Doob's $h$ and B\&S. We simulate 1000 Ornstein-Uhlenbeck 
            bridges. We vary the initial and endpoints.}
        
        \label{tb:Comparison_OU}
    \end{table}

\begin{figure}[H]
     \centering
      \subfloat[]{
        \includegraphics[width=0.48\textwidth,keepaspectratio]{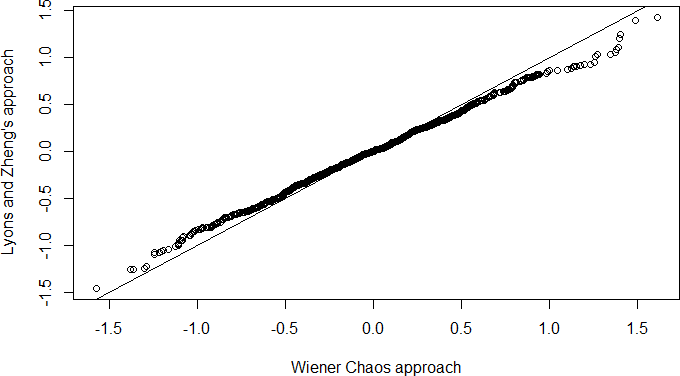}}
      \subfloat[]{
        \includegraphics[width=0.48\textwidth,keepaspectratio]{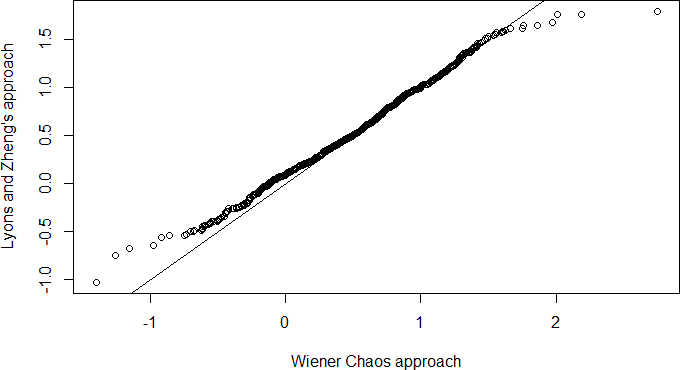}}
      \hfill
      \subfloat[]{
        \includegraphics[width=0.48\textwidth,keepaspectratio]{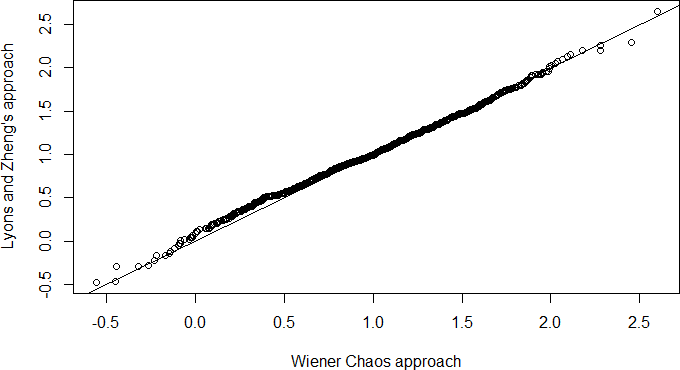}}
         \subfloat[]{
        \includegraphics[width=0.48\textwidth,keepaspectratio]{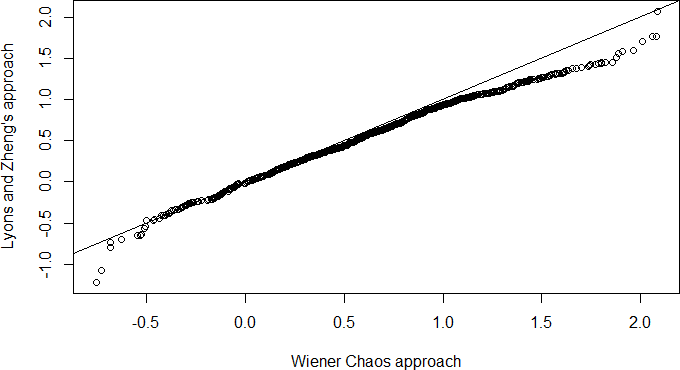}}
     \caption{Comparison of the empirical distributions of WCE and Doob's $h$-transform approach 
     at time 0.5. It is based on 1,000 realizations of Ornstein-Uhlenbeck bridges from  
     $(a)$ $(0,0)$ to $(1,0)$, $(b)$ from $(0,0)$ to $(1,1)$, $(c)$ from $(0,0)$ to 
     $(1,2)$ and $(d)$ from $(0,0.8)$ to $(1,0.5)$.}
     \label{f:q-q_ou2}
    \end{figure}

\begin{figure}[H]
     \centering
      \subfloat[]{
        \includegraphics[width=0.48\textwidth,keepaspectratio]{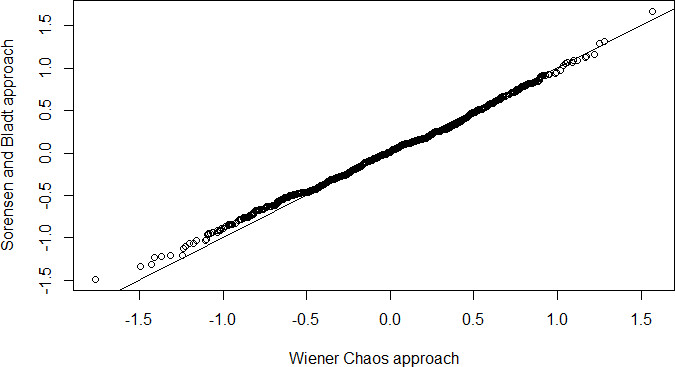}}
      \subfloat[]{
        \includegraphics[width=0.48\textwidth,keepaspectratio]{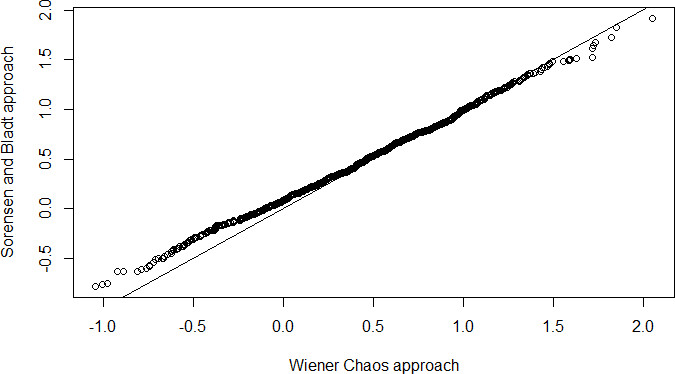}}
      \hfill
      \subfloat[]{
        \includegraphics[width=0.48\textwidth,keepaspectratio]{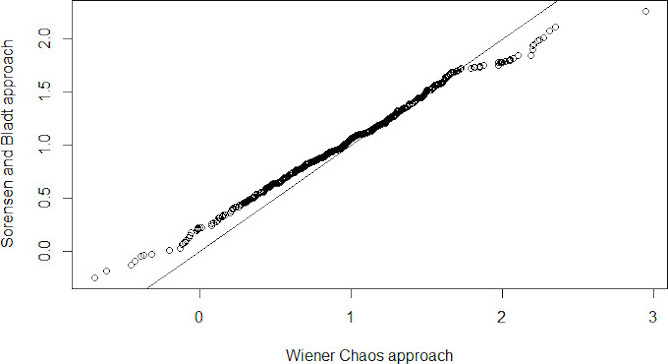}}
         \subfloat[]{
        \includegraphics[width=0.48\textwidth,keepaspectratio]{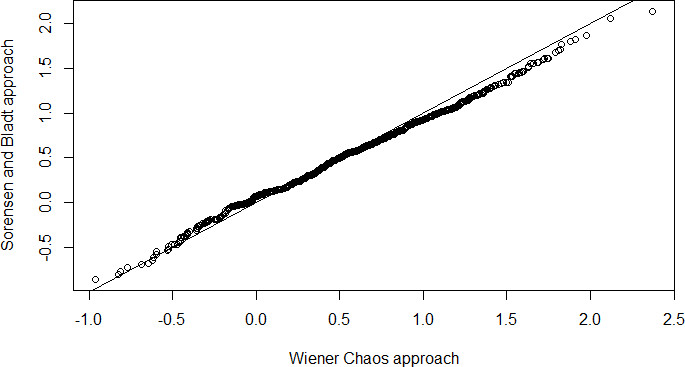}}
     \caption{Comparison of  the empirical distributions between WCE and S\o rensen and Bladt's approach at time 0.5. It is based on 1,000 realizations of Ornstein-Uhlenbeck bridges from  
     $(a)$ $(0,0)$ to $(1,0)$, $(b)$ from $(0,0)$ to $(1,1)$, $(c)$ from $(0,0)$ to 
     $(1,2)$ and $(d)$ from $(0,0.8)$ to $(1,0.5)$.}
     \label{f:q-q_ou}
    \end{figure}

    The results obtained from the WCE method and the approach by Bladt and S\o rensen appears to be comparatively weaker  (see \ref{other-method-app} for a brief description of the method). Although statistical evidence suggests that both methods have similar distributions in nearly all cases and simulations, it is noteworthy that as the size of the bridge increases, the p-value decreases. Specifically, in the case of the bridge from $(0,0) \to (1,2)$, the p-value falls below $0.05$, indicating a significant difference between the two distributions. This result could be explained due to some issues of the method of Bladt and S\o rensen. In some instances, the rejection rates can be extremely high, as is the case for the bridge from $(0,0) \to (1,2)$, due to the substantial gap. Furthermore, this specific bridge, constructed using this particular method, demonstrate a reduced amount of variation. Conversely, the WCE method consistently produces a bridge with a well-defined variance. This elucidates the reason why the calculated p-value is below the threshold of 0.05.
    We are also presenting the QQ-plots in Figures \ref{f:q-q_ou2} and \ref{f:q-q_ou} against  Doob's $h$-transform approach and  B\&S corresponding.

    As the quantiles from each QQ-plot align closely along a straight line, we can conclude that they came from the same distribution.
    
    \subsection{Geometric Brownian motion}
   
    Consider the following the geometric Brownian motion (GBM) that is the solution 
    of the SDE
    \begin{align}\label{SDE_GBM}
     dX(t)  & =     a X(t) dt + \sigma X(t) dB(t), \\
     X(0)  & =    x_0\nonumber.
    \end{align}
    We now describe the propagator for this SDE. For the case when $|m|=0$, we 
    denote by $X_0$ the solution of the propagator which solves the differential 
    equation
    \begin{align}\label{xm-equation0}
        X_0(t) 
            &= x_0 + a \int_0^t  X_0(s) ds. 
    \end{align}
    For the case $|m| \ge 1$ we have the vectors from \autoref{mi-values}. Then, $X_{m}^i$ satisfies
    \begin{align}\label{xm-equation2}
        X_{m}^i(t) 
            &= a \int_0^t  X_{m}^i(s) ds+ \sigma\sum_{j} 
            \sqrt{m_j^i} \int_0^t  e_j(s)\, X_{m^-(j)}^i(s)\, ds,
    \end{align}
    where $m^{i-}(j)$ is the $i$-th vector with the $j$-th entry minus $1$, see 
    the definition of  $m^{-}(j)$  in \eqref{m-i}.
    
    We observe that the exact solution of the equation \eqref{xm-equation0} is given by
    $$
    X_0(t)=x_0\exp(a t).
    $$
     Furthermore, the functions  $X_{m}^i(t)$, from expression \eqref{xm-equation2}, 
     with $|m|\ge 1$ will be solved numerically in increasing order of $|m|$. 
    
    We use the vectors $m^i$ described on \autoref{mi-values}, and we define 
    the propagator for each $m^i$ as 
    $$
    X_m^i(t)=\E\big(X(t)\xi_{m^i} \big)
    $$
     with $X(t)$ the solution to \eqref{SDE1}. We will consider the level of truncation as delineated in inequality \eqref{error-WCE}, to be $L = 0, 10, 100, 1000$, while maintaining a fixed value of $p=12$. 
     
\begin{figure}[H]
        \centering
        \includegraphics[width=1\textwidth, keepaspectratio]
        {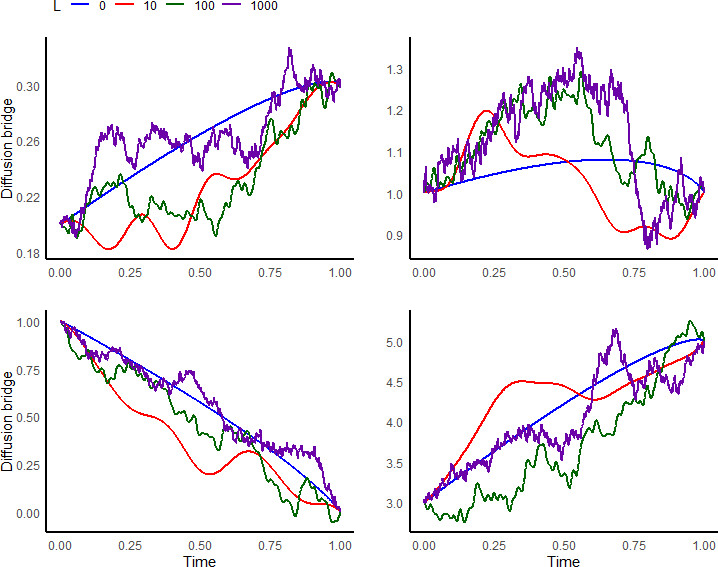}
        \caption{
            Geometric Brownian bridges with different truncation levels for the WCE,
            from a) (0,0.2) to (1, 0.3), b) (0, 1) to (1, 1), c) (0, 1) to (1, 0.1), and d) (0, 3) to (1, 5) with parameters $a = 0.2$, $\sigma = 0.3$.
        }
        \label{fig:Geometric_bridge}
    \end{figure}

     \Cref{fig:Geometric_bridge} illustrate fourth geometric Brownian bridges 
     using our method. Presenting the effect of increasing the truncation level are used in the Wiener approximation method. 
     The parameters for Geometric Brownian bridges are $a =0.5$, and 
     $\sigma =0.3$ with time discretization of 1000 point in the interval $[0,1]$.
    
\begin{figure}[H]
     \centering
      \subfloat[]{
        \includegraphics[width=0.5\textwidth,keepaspectratio]{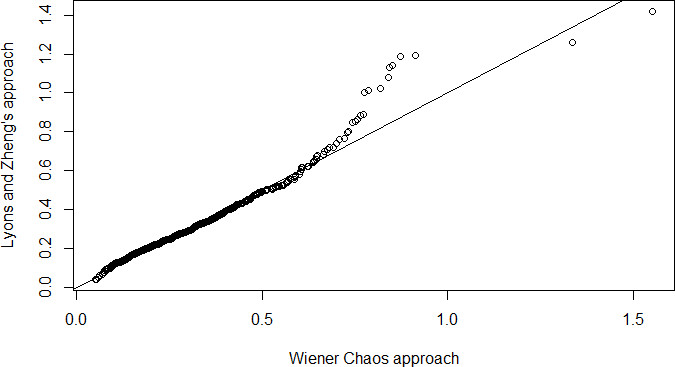}}
      \subfloat[]{
        \includegraphics[width=0.5\textwidth,keepaspectratio]{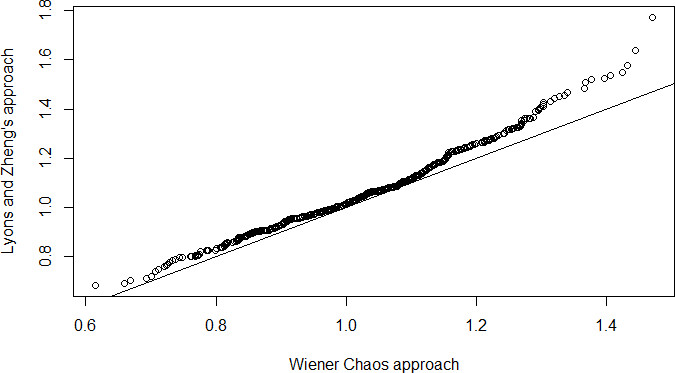}}
     \caption{Comparison of  the empirical distributions between WCE and  Doob's $h$-transform approach (see \ref{other-method-app}) 
     at time 0.5. It is based on 1,000 realizations of Geometric Brownian  bridges from  
     $(a)$ $(0,0.2)$ to $(1,0.3)$, $(b)$ from $(0,1)$ to $(1,1)$.}
     \label{f:q-q2}
\end{figure}

Given that the Geometric Brownian Motion does not possess a stationary distribution,  we are unable to apply the B\&S method. Consequently, we have restricted our comparison of results to Doob's $h$-transform approach exclusively (see \ref{other-method-app} for a brief description of the method). We present in Figure \ref{f:q-q2} a QQ-plot; from which we deduce that the results with both methods coincide. We also used the Kolmogorov-Smirnov test to compare the probability distributions from the Wiener chaos expansion method with those from  Doob's $h$-transform approach . Table \ref{tb:GBM} presents the results of this test. We assert, based on the statistical evidence presented in Table \ref{tb:GBM}, that the WCE method exhibits distributions comparable to those of the Doob's $h$-transform approach . Therefore, the WCE method constructs Geometric Brownian motion bridges. 

    \begin{table}[H]
    \centering
    \begin{tabular}{lr}
      \hline
     $(0,\eta) \to (T, \theta)$ & KS test  \\ 
     &  (p-value) \\ 
      \hline
      $(0,0.2) \to (1,0.3)$ & 0.0587\\
      $\,(0,1)\, \to \,(1,1)\,$ & 0.9270\\
     \hline
    \end{tabular}
        \caption{%
        This table contain the results of Kolmogorov-Smirnov test for the Wiener chaos expansion method and Doob's $h$-transform method. We perform 
             1000 Geometric Brownian bridges. We vary the initial and endpoints.
        }
        \label{tb:GBM}
    \end{table}

    \section{Application to two SDEs}\label{numerics2}
    \subsection{Multiplicative logistic SDE}
    
    Let us to define the stochastic multiplicative logistic process determined by the 
    following  SDE
    \begin{align}\label{SDE_Log}
     dX(t)  & = a X(t)\big(1-X(t)\big) dt + \sigma X(t) dB(t), \\
     X(0)  & =    x_0\nonumber.
    \end{align}
    Its propagator has the form given by the following expressions.
    
    We denote by $X_0$ the solution of the propagator for the case when $|m|=0$. Thus, 
    $X_0$ solves the differential equation
    \begin{align}\label{xm-log0}
        X_0(t) 
            &= x_0 + a \int_0^t  X_0(s)\big[1-X_0(s)\big] ds, 
    \end{align}
    which has as a solution
    \begin{align}\label{xm-log0-sol}
        X_0(t) 
            &=\frac{x_0}{1+\tfrac{1-x_0}{x_0} e^{-a t}}. 
    \end{align}
    
    We now consider the case where $|m| = 1$, that, when applied to the vectors from
    \autoref{mi-values}. Thus, $X_{m}^i$ satisfies the following equation.
    
    \begin{align}\label{xm-log1}
        X_{m}^i(t) 
            &= a \int_0^t  X_{m}^i(s) \big[1- X_{m}^i(s) \big] ds+ \sigma\sum_{j} 
            \int_0^t e_j(s)\, X_{m^{i-}(j)}^i(s)\ ds,
    \end{align}
    
    where $m^{i-}(j)$ is the $i$-th vector with the $j$-th entry minus $1$ (see the 
    definition of $m^{-}(j)$ in \eqref{m-i}).
    
    Finally, for the case with $|m| > 1$ for the vectors from \autoref{mi-values}, 
    $X_{m}^i$ satisfies
    \begin{align}\label{xm-log2}
        X_{m}^i(t) 
            &= a \int_0^t  X_{m}^i(s) \big[1- X_{m}^i(s) \big]  ds 
                + \sigma\sum_{j} \sqrt{m_j^i} 
            \int_0^t e_j(s)\, X_{m^-(j)}^i(s)\, ds,
    \end{align}
    where $m_j^i$ denotes the $j$-th entry of the $i$-th vector.
    
    We observe, as in the precedent example, that the functions  $X_{m}^i(t)$ with 
    $|m|\ge 1$ will be solved numerically using expressions \eqref{xm-log1} and 
    \eqref{xm-log2} in increasing order. We will consider the level of truncation to be $L = 0, 10, 100, 1000$, with a fixed value of $p=12$. 
    
\begin{figure}[H]
        \centering
        \includegraphics[width=1\textwidth, keepaspectratio]
        {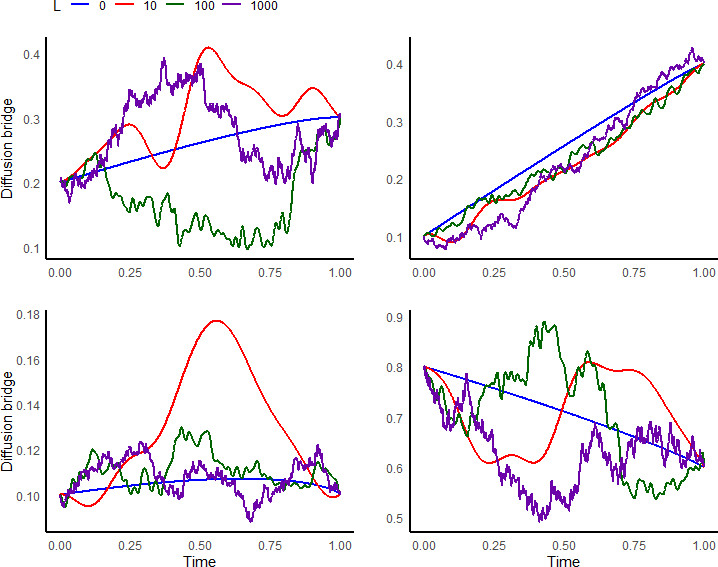}
        \caption{
            Multiplicative logistic bridges with different truncation levels for the WCE from a) (0,0.2) to (1, 0.3), b) (0, 0.1) to (1, 0.4), 
            c) (0, 0.1) to (1, 0.1), and d) (0, 0.8) to (1, 0.6). The parameters of the SDE are
            $a = 0.2$, $\sigma = 0.7$.
        }
        \label{fig:MultLogistic_bridges}
    \end{figure}

    \Cref{fig:MultLogistic_bridges} presents four cases of Multiplicative 
    Logistic Bridges using Wiener approximation method. We show the effect 
    of increasing the number of Brownian motions. We take $dt = 0.001$ in 
    the interval $[0,1]$ and parameters $a =0.5$, and $\sigma =0.7$.
     
\subsection{Stochastic protein kinetic differential equation}
    
    Finally, we consider the following protein kinetic SDE in Stratonovich sense
    \begin{align}\label{SDE_Kin_stra}
     dX(t)  & = \Big[ 1-X(t) + \lambda X(t) \big(1-X(t)\big)\Big] dt + \sigma X(t)\big(1-X(t)\big) 
        \circ dB(t), \\
     X(0)  & =    x_0\nonumber.
    \end{align}    
    This equation describes the kinetics of the proportion $X$ of one of two 
    possible forms of certain proteins where $\lambda$ is the interaction coefficient 
    of the two proteins (see for instance \cite{wang2016effective}).
    
    Rewriting \eqref{SDE_Kin_stra} in It\^o sense we obtain
    
    \begin{align}\label{SDE_Kin}
     dX(t)  & =     \Big[ 1-X(t) + \lambda X(t) \big(1-X(t)\big)  
     + \sigma X(t)\big(1-X(t)\big) \big(1-2X(t)\big) \Big] dt\nonumber\\ 
     &\qquad + \sigma X(t)\big(1-X(t)\big) dB(t), \\
     X(0)  & =    x_0\nonumber.
    \end{align}
    
    We now write the propagator for this SDE. We denote by $X_0$ the propagator for 
    the case $|m|=0$. Thus, 
    \begin{align}\label{kin_prop0}
    X_0(t)= (x_0-1) + \int_0^t    \Big[  (\lambda+\sigma-1) X_0(s) 
    - (\lambda+3\sigma) \big(X_0(s))^2  + 2 \sigma  \big(X_0(s))^3 \big) 
    \Big] ds.
    \end{align}
    
    We now consider the case where $|m| = 1$, that, when applied to the vectors from
    \autoref{mi-values}. Thus, $X_{m}^i$ satisfies 
    the following equation.
    
    \begin{align}\label{kin_prop1}
    X_m^i(t)&= \int_0^t    \Big[  (\lambda+\sigma-1) X_m^i(s) -(\lambda+3\sigma) 
    \big(X_m^i(s))^2  + 2 \sigma  \big(X^i(s))^3 \big) \Big] ds \nonumber\\
    &\qquad + \sigma \sum_{j} \int_0^t  e_j(s) \Big[  X_{m-(j)}^i(s) 
    \big( 1- X_{m-(j)}^i(s) \big) \Big] ds,
    \end{align}
    
    where, as before, $m^{i-}(j)$ is the $i$-th vector with the $j$-th entry minus $1$ 
    (see the definition of  $m^{-}(j)$  in \eqref{m-i}).
    
    Finally, for the case with $|m| > 1$ we have the vectors from \autoref{mi-values} 
    and then $X_{m}^i$ satisfies
    \begin{align}\label{kin-prop2}
        X_m^i(t)&= \int_0^t    \Big[  (\lambda+\sigma-1) X_m^i(s) 
        -(\lambda+3\sigma) \big(X_m^i(s))^2  
        + 2 \sigma \big(X^i(s))^3 \big) \Big] ds \nonumber\\
        &\qquad + \sigma \sum_{j} \sqrt{m_j^i} \int_0^t  e_j(s) 
        \Big[  X_{m-(j)}^i(s) \big( 1- X_{m-(j)}^i(s) \big) \Big] ds,
    \end{align}
    where $m_j^i$ denotes the $j$-th entry of the $i$-th vector.
    
As before, we will consider the level of truncation to be $L = 0, 10, 100, 1000$, with a fixed value of $p=12$. 
    
\begin{figure}[H]
        \centering
        \includegraphics[width=1\textwidth, keepaspectratio]
        {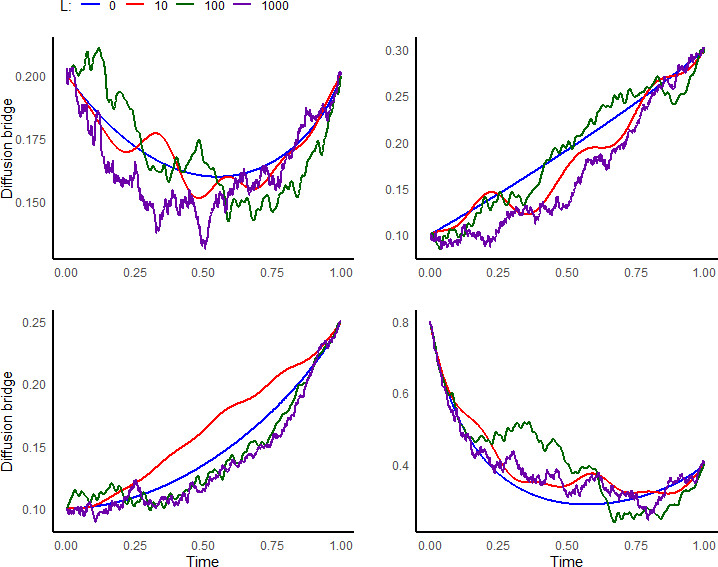}
        \caption{
            Protein Kinetic bridges with different truncation leves of the WCE from a) (0,0.2) to (1, 0.2), b) (0, 0.1) to (1, 0.3), 
            c) (0, 0.1) to (1, 0.25), and d) (0, 0.8) to (1, 0.4). The parameters of the SDE are
            $a = 0.2$, $\sigma = 0.8$.
        }
        \label{fig:ProteinKinetic_bridges1}
    \end{figure}

    In Figure \ref{fig:ProteinKinetic_bridges1}, we show Protein Kinetic Bridges 
    using the Wiener chaos approximation method. Illustrating the effect of increasing 
    the truncation level $L$ in the Wiener chaos expansion. In this figure, 
    we consider the parameters  $a =0.5$, and $\sigma =0.8$, with 1000 points in the 
    interval $[0,1]$.
    
    \section{Concluding remarks}\label{conclusions}
In this paper, we propose a numerical method to approximate diffusion bridges efficiently. This indicates that regardless of the size of the gap, whether in terms of time or state space, our methodology consistently establishes a bridge with a uniform computational cost. The proposed numerical method is based on a truncation of the WCE for a particular SDE. In order to implement the proposed simulation method, it is necessary to solve a set of coupled ODEs.

    We have validated the method with the OU bridges. Indeed, in the simulation study, we showed that the distribution of the OU bridges, obtained with the WCE method, is the same distribution as that from exact simulations. Moreover, we found numerical evidence that the execution time is linear by increasing the number of Brownian motions in the WCE expansion.
    
   We have also shown that our method works well in the cases when the coupled equation system has to be solved numerically.
    
     The methodology delineated in this paper has been successfully implemented for SDEs in one dimension; however, it possess the potential for extension to higher dimensions in a relatively straightforward manner. The challenge associated with this extension appears to be primarily concerned with the solution, likely in a numerical sense, of the system of ODEs that possesses the same dimension as the corresponding system of SDEs. Upon obtaining the numerical solution, we will be positioned to construct multidimensional diffusion bridges. This will be the focus of a future project.
    
    We will also examine diffusion bridges driven by fractional Brownian motion, which is more complex than those driven by classical Brownian motion. However, it appears possible to conduct simulations of bridges using the ideas presented here.
    
Our proposed method has an additional advantage in data augmentation; allow us to explain further. As we mentioned before, one advantage of our method, over previous bridge simulation methods, is that it works very well and with the same speed of execution for diffusion bridges; this is true even in long intervals. In addition, assume the situation when the true value of the parameters in the SDE is unknown and it is necessary to use an estimator to perform simulations; this scenario is common when we use an EM-algorithm or an MCMC method to find the estimators. Thus, it is known that if the estimator is not close enough to the true value, then other simulation methods could have high rejection rates when constructing diffusion bridges. In our method, that is not the case; indeed, our method allows us to construct a bridge even if the parameter values (that we used in the bridge simulation) differ significantly from the true values. Hence, the methodology employed in this study provides a crucial advantage in data augmentation.  In future research, we will study the effectiveness of our method for data augmentation.

    \appendix
    
    \section{Concrete choices of vectors in $\mathcal{J}^{p,L}$}\label{app-B}
    In this section, we present explicitly the vectors we use for the  
    numerical examples. We define the set
    $\mathcal{J}^{p,L}$ as is shown in Table \ref{mi-values}. 
    
    We will use the notation $m^i$ to denote the  $i$-th vector where 
    $|m^i|\leq p$ (maximum order of the Hermite polynomial used), and 
    with a number of entries $L$.
    
    \begin{table}[H]
    \centering
    \begin{tabular}{rrrrrrrrr|r}
      \hline
     & $m^i(1)$ &$m^i(2)$ & $m^i(3)$ & $m^i(4)$ & $m^i(5)$ &$m^i(6)$ & $\cdots $ & $m^i(L)$ & $|m^i|$  \\ 
      \hline
      $m^{0}$ & 0 & 0 &  0 & 0 & 0 & 0 & $\cdots$ &0
      &  0\\ 
      $m^1$ & 1 &  0 & 0 &0& 0 & 0 & $\cdots$ &0 & 1 \\ 
      $m^2$ & 0 & 1 & 0 & 0 & 0 &  0& $\cdots$ &0  & 1\\ 
      $m^3$ & 0 & 0 & 1 & 0 & 0 &  0& $\cdots$ &0  & 1\\ 
      $m^4$ & 0 & 0 & 0 & 1 & 0&  0& $\cdots$ &0  & 1\\ 
      $m^5$ & 0 & 0 & 0 & 0 & 1 & 0&  $\cdots$ &0 & 1\\ 
      $m^6$ & 0 & 0 & 0 & 0 & 0 & 1 &  $\cdots$ &0  & 1\\ 
      $\vdots$ &  $\vdots$ &  $\vdots$ &  $\vdots$ &   $\vdots$ 
      &  $\vdots$ &  $\vdots$ & $\vdots$ &  $\ddots$ & 1  \\
       $m^{j}$ & 0 & 0 & 0 & 0 & 0 & 0&  $\cdots$ &1 & 1 \\ 
         $m^{j+1}$ & 1 & 1 & 0 & 0 & 0 & 0&  $\cdots$ &0 & 2\\ 
       $m^{j +2}$ & 1 & 0 & 1 & 0 & 0 & 0& $\cdots$ &0 & 2 \\ 
      $m^{j +3}$ & 0 & 1 & 1 & 0 & 0 & 0& $\cdots$ &0  & 2\\ 
      $m^{j +4}$ & 2 & 0 & 0 & 0 & 0& 0&  $\cdots$ &0 & 2 \\ 
      $m^{j +5}$ & 0 & 2 & 0 & 0 & 0 & 0&  $\cdots$ &0 & 2 \\ 
      $m^{j +6}$ & 0 & 0 & 2 & 0 & 0& 0 &  $\cdots$ &0 & 2\\ 
      $m^{j +7}$ & 1 & 2 & 0 & 0 & 0 & 0&  $\cdots$ &0 & 3\\ 
      $m^{j +8}$ & 2 & 1 & 0 & 0 & 0& 0 &  $\cdots$ &0 & 3\\ 
      $m^{j +9}$ & 3 & 0 & 0 & 0 & 0& 0&  $\cdots$ &0 & 3\\ 
      $m^{j +10}$ & 0 & 3 & 0 & 0 & 0& 0&  $\cdots$ &0 & 3\\ 
      $m^{j +11}$ & 4 & 0 & 0 & 0 & 0 & 0 &  $\cdots$ &0 & 4\\ 
      $m^{j +12}$ & 0 & 4 & 0 & 0 & 0 & 0&  $\cdots$ &0 & 4\\
      $m^{j +13}$ & 5 & 0 & 0 & 0 & 0 & 0&  $\cdots$ &0 & 5\\
      $m^{j +14}$ & 0 & 5 & 0 & 0 & 0 & 0&  $\cdots$ &0 & 5\\
      $m^{j +15}$ & 6 & 0 & 0 & 0 & 0 & 0&  $\cdots$ &0 & 6\\
       $m^{j +16}$ & 0 & 6 & 0 & 0 & 0 & 0&  $\cdots$ &0 & 6\\
       $m^{j +17}$ & 7 & 0 & 0 & 0 & 0 & 0&  $\cdots$ &0 & 7\\
       $m^{j +18}$ & 0 & 7 & 0 & 0 & 0 & 0&  $\cdots$ &0 & 7\\
       $m^{j +19}$ & 8 & 0 & 0 & 0 & 0 & 0&  $\cdots$ &0 & 8\\
       $m^{j +20}$ & 0 & 8 & 0 & 0 & 0 & 0&  $\cdots$ &0 & 8\\
         $m^{j +21}$ & 9 & 0 & 0 & 0 & 0 & 0&  $\cdots$ &0 & 9\\
       $m^{j +22}$ & 0 & 9 & 0 & 0 & 0 & 0&  $\cdots$ &0 & 9\\
       $m^{j +23}$ & 10 & 0 & 0 & 0 & 0 & 0&  $\cdots$ &0 & 10\\
       $m^{j +24}$ & 0 & 10 & 0 & 0 & 0 & 0&  $\cdots$ &0 & 10\\
     \hline
    \end{tabular}
        \caption{%
            Table with the vectors $m^i\in \mathcal{J}^{p,L}$ we have used to perform the numerical experiments. Here 
            $m^i(k)$ denotes the $k$-th entry of the $i$-th vector, where $|m^i|\leq p$. 
            $j$ is taken large enough to ensure a good approximation 
            to the bridge, see \eqref{trunc_WCE_X}.
            }
        \label{mi-values}
    \end{table}

\section{Approachs for construction of diffusion bridges}\label{other-method-app}

In this section, we provide two different numerical construction of diffusion bridges that we use to compare the method proposed in this work.

\subsection{Doob's $h$-transform approach}

In \cite{delyon2006simulation}, the authors use a result from an article written by Lyons and Zheng (see \cite{lyons1990conditional}), to propose some algorithms for the simulation of diffusion bridges. Indeed, Lyons and Zheng show that the distribution of the diffusion bridge given by \eqref{SDE1}-\eqref{SDE_conditions} is the same as the diffusion 
\begin{align}
dz_t &= b(z_t) dt + \sigma(z_t) dB(t)\\
z_0&= \eta,
\end{align}
where 
$$
b(z)= f(z) + \sigma(z) \sigma^T(z)\, \nabla_z\big(\log\, p(t,z;T,\theta) \big) 
$$
and $p(t,z;T,\theta)$ is the transitional probability density of the diffusion bridge, evaluated at
terminal time $T$. It is recognized that only a limited number of stochastic differential equations  possess a known transition density. Consequently, employing this method for simulations is not deemed appropriate. Nevertheless, this approach is employed to compare the method presented in this manuscript when dealing with the OU process as well as with the geometric Brownian motion, as it is known these possesses an explicit transition density. The method used for this bridge simulation will be referred to as the {\it Doob's $h$-transform approach}.

\subsection{Bladt and S\o rensen's approach}

The method proposed in \cite{bladtSoren} and \cite{bladtSorenc} is based on the following simple construction: one diffusion process $\bm X^1$ solution of \eqref{SDE1}, is simulated using as the starting point $\eta$ ($X_t^1=\eta$), while another independent solution to \eqref{SDE1}, $\bm X^2$, is started from the point $\theta$ ($X_t^2=\theta$). The time of the second diffusion is then reversed to obtain the simulated process $X^r_t=X^2_{T-1}$. Suppose there is a time point $\tau\in[0,T]$ at which the paths satisfies $X^1_{\tau} = X^r_{\tau}$. Then, they propose as a {\it diffusion bridge} the process that is equal to $X^1_t$ for $t\in [0,\tau]$ and $X^r_t$ for $t\in[\tau,T]$: that is a process that starts at $\eta$ and ends at $\theta$.  It is important to remark that the probability that $\bm X^1$ and $\bm X^r$ do not meet in the interval $[0,T]$ is significant, as long as $T$ or $\vert\theta-\eta\vert$ are not very small. Consequently, this method may not be efficient.

    \section*{Software and data availability}
        The data, Julia, and R source code are available at
    
        \href{https://github.com/AdrianSalcedo/DiffusionBridgesWienerExpantion.git}
        {https://github.com/AdrianSalcedo/DiffusionBridgesWienerExpantion.git}

   \section*{Acknowledgements} 
F. Baltazar-Larios and F. Delgado-Vences have been supported by UNAM-DGAPA-PAPIIT-IN102224 (Mexico). G.A. Salcedo-Varela  is supported by "Estancias Postdoctorales por 
M\'exico-Iniciales" CVU: 711305 (CONAHCYT).\\

\bibliographystyle{elsarticle-num}
\bibliography{references}
\end{document}